\documentclass[11pt,oneside]{article}

\usepackage{amscd}
\usepackage{amsfonts}
\usepackage[centertags]{amsmath}
\usepackage{amssymb}
\usepackage[english]{babel}
\usepackage[all]{xy}
\usepackage{theorem}

\newcommand{\al}{\ensuremath{\alpha}}
\newcommand{\be}{\ensuremath{\beta}}

\newcommand{\T}{\ensuremath{\theta}}

\newcommand{\map}[3]{\mbox{${#1}\colon{#2}\to{#3}$}}
\newcommand{\dismap}[5]{
\[
\begin{array}{rcll}
#1: & #2 & \longrightarrow & #3\\
 & #4 & \longmapsto & #5
\end{array}
\]
}
\newcommand{\di}{\displaystyle}

\newcommand{\R}{\ensuremath{\mathbb R}}
\newcommand{\C}{\ensuremath{\mathbb C}}

\newcommand{\Ha}{\ensuremath{\mathbb H}}

\newcommand{\Cp}[1]{\ensuremath{{\C\mathbb P}^{#1}}}

\newcommand{\ug}{\stackrel{\rm def}{=}}

\newcommand{\iso}{\ensuremath{\simeq}}

\newcommand{\myref}[1]{(\ref{#1})}

\newcommand{\EndDim}{\ensuremath{\nopagebreak\hfill\blacksquare}}
\newcommand{\EndRem}{\ensuremath{\nopagebreak\hfill\square}}
\newcommand{\EndDef}{\ensuremath{\nopagebreak\hfill\square}}

\newcommand{\pa}[1]{\ensuremath{\partial_{#1}}}
\newcommand{\copa}[1]{\ensuremath{dx_{#1}}}

\newcommand{\st}{\ensuremath{\text{ such that }}}

\renewcommand{\setminus}{-}

\newtheorem{lem}{Lemma}[section]

\newtheorem{teo}[lem]{Theorem}
\newtheorem{cor}[lem]{Corollary}
\newtheorem{pro}[lem]{Proposition}
{\theorembodyfont{\rmfamily} \newtheorem{oss}[lem]{Remark}
                             \newtheorem{defi}[lem]{Definition}
                             
                            }

\newenvironment{D}[1][]{{\nopagebreak\em Proof#1: }}{\EndDim}
\newenvironment{Dsketch}[1][]{{\nopagebreak\em Sketch of proof#1: }}{\EndDim}
\newenvironment{acknowledgements}{{\em Acknowledgements: }}{}

\author{ Maurizio Parton,\\
 Dipartimento di Matematica ``Leonida Tonelli'',\\
 via Filippo Buonarroti 2, I--56127 Pisa,\\
 Italy. E-mail:  {\tt parton@dm.unipi.it}
}
\date{}

\begin{document}

\title{Explicit parallelizations on products of spheres}


\maketitle


\section{Introduction}

It is a classical result in Algebraic Topology that spheres $S^n$
are
parallelizable if and only if their dimension is $n=1,3,7$.
As for the products of two or more spheres
one has instead the following
result (\cite{KerCIG}):

\begin{teo}[Kervaire]\label{teokervaire}
The manifold $S^{n_1}\times\dots\times S^{n_r}$, $r\ge 2$, is
parallelizable if and only if at least one of the $n_i$ is odd.
\end{teo}

Kervaire's proof does not provide an explicit parallelization on
products of spheres.
The only reference the
author knows to provide explicit parallelizations
is \cite{BruPEP}, that considers the cases when one of the spheres is of
dimension $1,3,5,7$,
and uses some specific arguments of these low dimensions. In \cite{BruPEP} the
general case is left as an open problem.

%

The aim of this paper is to write an explicit orthonormal
parallelization for all
parallelizable products of spheres, using an explicit
isomorphism with a trivial vector bundle obtained
following a hint of \cite{HirDTo}.

%
%

A description of some $G$-structures on $S^m\times S^n$
associated to this parallelization is given in \cite{ParPhT}.
%
%
%

\section{An explicit parallelization $\mathcal{B}$ on $S^m\times S^1$}\label{genm1}

Denote by $x=(x_i)$
the coordinates on $\R^{m+1}$, and let $S^m\subset\R^{m+1}$ be
given by
\[
\begin{split}
S^m&\ug\{x=(x_1,\dots,x_{m+1})\in\R^{m+1}\st
|x|^2=x_1^2+\dots+x_{m+1}^2=1\}.\\
\end{split}
\]
The orthogonal projection of the standard coordinate
frame $\{\pa{x_i}\}_{i=1,\dots,m+1}$ to the sphere plays an
important role in the game, and deserves its own definition:
\begin{defi}
The {\em $\imath^{\text{th}}$ meridian vector field} $M_i$
on $S^m$ is 
\[
\begin{split}
M_i&\ug\text{orthogonal
projection of $\pa{x_i}$ on $S^m$}\qquad i=1,\dots,m+1.\\
\end{split}
\]\EndDef
\end{defi}
Let $M$ be
the normal versor field of $S^m\subset \R^{m+1}$, that is,
\[
M\ug\sum_{i=1}^{m+1}x_i\pa{x_i}.
\]
Since
\[
\langle \pa{x_i},M\rangle=x_i \qquad i=1,\dots,m+1,
\]
one obtains the following expression for $M_i$:
\begin{equation}\label{fundeq}
M_i=\pa{x_i}-x_i M\qquad i=1,\dots,m+1,
\end{equation}
and thus
\begin{equation}\label{fundeqscal}
\langle M_i,M_j\rangle=\delta_{ij}-x_i x_j\qquad i,j=1,\dots,m+1.
\end{equation}



Let $\Gamma$ be the cyclic infinite group of transformations
of $\R^{m+1}\setminus 0$ generated by the map $x\mapsto e^{2\pi}x$.
Denote by $H$ the
corresponding diagonal real Hopf manifold, that is,
the quotient manifold $(\R^{m+1}\setminus 0)/\Gamma$:
$H$ turns out to be diffeomorphic to $S^m\times
S^1$ by means of the map
induced by the projection $p$:
\[
\begin{array}{cll}\label{projm1}
\R^{m+1}\setminus 0 & \stackrel{p}{\longrightarrow} & S^m\times S^1\\
 x & \longmapsto &
 (x/|x|,\log|x| \mod{2\pi}).
\end{array}
\]



The standard coordinate frame $\{\pa{x_i}\}_{i=1,\dots,m+1}$ on $\R^{m+1}\setminus
0$ becomes $\Gamma$-equivariant when multiplied by the function $|x|$,
whence it defines a parallelization on $S^m\times S^1$.
This proves the following proposition\dots
\begin{pro}
$S^m\times S^1$ is parallelizable.
\end{pro}
\dots and enables us to give the following definition:
\begin{defi}\label{Bdef}
Define
${\cal B}=\{b_i\}_{i=1,\dots,m+1}$
as the frame  on $S^m\times S^1$ induced
by the $\Gamma$-equivariant frame $\{|x|\pa{x_i}\}_{i=1,\dots,m+1}$ on the
universal covering $\R^{m+1}$ of $S^m\times S^1$ by means of $p$:
\[
b_i\ug
p_*(|x|\pa{x_i}(x))
\qquad{i=1,\dots,m+1}.
\]\EndDef
\end{defi}



The following theorem explicitly describes the frame $\mathcal{B}$:
\begin{teo}\label{progeo}
Let $M_i$ be the $i^{\text{th}}$ meridian vector field
on $S^m\subset\R^{m+1}$.
Then
\begin{equation}\label{progeofor}
b_i=M_i+x_i\pa{\T} \qquad i=1,\dots,m+1.
\end{equation}
\end{teo}

\begin{D}
Look at $S^m\times S^1$ as a Riemannian submanifold
of $R^{m+1}\times S^1$, and in particular look at $T(S^m\times
S^1)=TS^m\times TS^1$ as a Riemannian subbundle of
$TR^{m+1}_{|_{S^m}}\times
TS^1$; this last is a trivial vector bundle and an orthonormal frame is
$\{\pa{x_1},\dots,\pa{x_{m+1}},\pa{\T}\}$. A computation then
shows  that
\begin{equation*}
\begin{split}
p_*=\di\frac{1}{|x|}\left((dx_1-x_1\omega)\otimes\pa{x_1}+\dots+
(dx_{m+1}-x_{m+1}\omega)\otimes\pa{x_{m+1}}+|x|\omega\otimes \pa{\T}\right),
\end{split}
\end{equation*}
where  $\omega$ is
the 1-form
given  by
\[
\omega\ug -d\left(\di\frac{1}{|x|}\right)=
\di\frac{1}{|x|^2}\left(x_1dx_1+\dots+x_{m+1}dx_{m+1}\right).
\]
Whence, the frame $\cal B$
in the point $p(x)=(x/|x|,\log|x| \mod{2\pi})$
is given by
\begin{equation*}
\begin{split}
\di\frac{1}{|x|^2}\left(-x_1x_i\pa{x_1}+\dots+(|x|^2-x_i^2)\pa{x_i}
+\dots-
x_{m+1}x_i\pa{x_{m+1}}+|x|x_i\pa{\T}\right)
\end{split}\qquad i=1,\dots,m+1,
\end{equation*}
that is,
the frame $\cal B$
in the point $(x,\T)\in S^m\times S^1$ is given by
\begin{equation}\label{bruni1}
\begin{split}
b_i&=\left(-x_1x_i\pa{x_1}+\dots+(1-x_i^2)\pa{x_i}+\dots-
x_{m+1}x_i\pa{x_{m+1}}+x_i\pa{\T}\right)\\
&=\pa{x_i}-x_i(x_1\pa{x_1}+\dots+
x_{m+1}\pa{x_{m+1}})+x_i\pa{\T}\stackrel{\myref{fundeq}}{=} M_i+x_i\pa{\T}.
\end{split}\qquad i=1,\dots,m+1.
\end{equation}
\end{D}


\begin{oss}
The notion of meridian vector field was given in
\cite{BruPEP}: it was used to describe a
parallelization on any product of  a
sphere by a parallelizable manifold. In this context, theorem
\ref{progeo} shows that the frame $\mathcal{B}$ given by definition \ref{Bdef}
coincide with that of \cite{BruPEP}.\EndRem
\end{oss}

\begin{oss}\label{remor}
The frame $\cal B$ is orthonormal with respect to the product
metric on $S^m\times S^1$ (use theorem \ref{progeo} and formula \myref{fundeqscal}).\EndRem
\end{oss}


The well-known bracket formula 
\begin{equation}\label{kobano}
[fX,gY]=fg[X,Y]+f(Xg)Y-g(Yf)X
\end{equation}
gives the brackets of $\cal B$:
\begin{equation}\label{bracket1gen}
[b_i,b_j]=x_ib_j-x_jb_i\qquad i,j=1,\dots,m+1.
\end{equation}

By means of theorem \ref{progeo} and remark \ref{remor},
one obtains the coframe ${\cal B}^*\ug\{b^i\}_{i=1,\dots,m+1}$ dual to  $\cal B$ on
$S^m\times S^1$:
\begin{equation}\label{cobruni1}
\begin{split}
b^i&= \copa{i}+x_i d\T
\end{split}\qquad i=1,\dots,m+1.
\end{equation}

\begin{oss}
Since
\[
b_i=p_*(|x|\pa{x_i})\qquad i=1,\dots,m+1,
\]
the coframe ${\cal B}^*$ can be also described as the quotient of
the
$\Gamma$-invariant coframe on $\R^{m+1}\setminus 0$ given by
\[
\{|x|^{-1}\copa{i}\}_{i=1,\dots,m+1}.
\]\EndRem
\end{oss}

A straightforward computation gives the structure equations for $\cal
B$:
\begin{equation}
db^i=\copa{i}\wedge d\T\stackrel{\myref{cobruni1}}{=}b^i\wedge{d\T}\qquad i=1,\dots,m+1,
\end{equation}
where the 1-form $d\T$ is related to ${\cal B}^*$ by
\[
d\T=\sum_{i=1}^{m+1}x_ib^i.
\]

The following lemma is trivial to prove, but will be useful:
\begin{lem}\label{lemgen}
For each permutation $\pi$ of $\{1,\dots,m+1\}$, the automorphism
of $\R^{m+1}\setminus 0$ given by $(x_1,\dots,x_{m+1})\mapsto
(x_{\pi(1)},\dots,x_{\pi({m+1})})$ is $\Gamma$-equivariant. The induced
diffeomorphism is
\dismap{f_\pi}{S^m\times S^1}{S^m\times
S^1}{(x_1,\dots,x_{m+1},\T)}{(x_{\pi(1)},\dots,x_{\pi({m+1})},\T),}
and $df_\pi(b_{\pi(i)})=b_i$.
\end{lem}

\section{The Hopf fibration $S^3\rightarrow S^2$ extends $\mathcal{B}$ to
$S^m\times S^3$}\label{genm3}

Denote by $y=(y_j)$
the coordinates on $\R^{4}$, and let $S^3\subset\R^{4}$ be
given by
\[
\begin{split}
S^3&\ug\{y=(y_1,\dots,y_{4})\in\R^{4}\st
|y|^2=y_1^2+\dots+y_{4}^2=1\}.\\
\end{split}
\]

Let $T=T_1,T_2,T_3$ be the vector fields on $S^3$ given
by multiplication by $i,j,k\in\Ha=\R^4$ respectively, that is,
\begin{equation}
\begin{split}
T=T_1&=-y_2\pa{y_1}+y_1\pa{y_2}-y_4\pa{y_3}+y_3\pa{y_4},\\
T_2&=-y_3\pa{y_1}+y_4\pa{y_2}+y_1\pa{y_3}-y_2\pa{y_4},\\
T_3&=-y_4\pa{y_1}-y_3\pa{y_2}+y_2\pa{y_3}+y_1\pa{y_4}.
\end{split}
\end{equation}

It is well-known that $S^3$ can be foliated in $S^1$'s, by means
of the Hopf fibration $S^3{\rightarrow}S^2$. Whence,
one has a foliation of $S^m\times S^3$ in $S^m\times S^1$'s, and
section \ref{genm1} gives $m+1$ vector fields tangent to the
leaves: they can be completed to a parallelization of $S^m\times S^3$
by means of a suitable parallelization of $S^3$, as it is now going
to be shown in the following proposition:

\begin{pro}\label{parm3}
$S^m\times S^3$ is parallelizable.
\end{pro}

\begin{D}
One would like to use definition \ref{Bdef}
to define $m+1$ vector fields on $S^m\times S^3$. The problem
is that there is not a canonical identification of the fiber of
the Hopf fibration $S^3{\rightarrow}S^2$ with $S^1$,
whence one has not a canonical angular coordinate on fibers.
But formula \myref{progeofor}
 of theorem \ref{progeo}
only needs a unitary and tangent to fibers vector field on $S^3$
to be used: this is just what $T$ is.
Whence, define
${\cal B}\ug\{b_i\}_{i=1,\dots,m+3}$ by
\begin{equation}\label{parsms3}
\begin{split}
b_i&\ug M_i+x_i T\qquad i=1,\dots,m+1,\\
b_{m+j}&\ug T_j\qquad j=2,3,
\end{split}
\end{equation}
where $M_i$ is the $\imath^{\text{th}}$ meridian vector field on
$S^m$,
to obtain the wished frame on $S^m\times S^3$.
\end{D}

\begin{oss}\label{remor3}
The frame $\cal B$ is orthonormal with respect to the product
metric on $S^m\times S^3$ (use formulas \myref{parsms3} and formula \myref{fundeqscal}).\EndRem
\end{oss}

The same argument used in section \ref{genm1} gives the brackets of $\cal
B$:
\begin{equation}\label{bracket3gen}
\begin{split}
[b_i,b_j]&=x_ib_j-x_jb_i\qquad i,j=1,\dots,m+1,\\
[b_i,b_{m+2}]&=-2x_ib_{m+3}\qquad i=1,\dots,m+1,\\
[b_i,b_{m+3}]&=2x_ib_{m+2}\qquad i=1,\dots,m+1,\\
[b_{m+2},b_{m+3}]&=-2T=-2\sum_{i=1}^{m+1}x_ib_i.
\end{split}
\end{equation}

Let $\tau=\tau_1,\tau_2,\tau_3$ be the 1-forms on $S^m\times S^3$ dual to
$T=T_1,T_2,T_3$ respectively.
The coframe ${\cal B}^*\ug\{b^i\}_{i=1,\dots,m+3}$
is given by
\begin{equation}\label{cobruni3}
\begin{split}
b^i&= x_i \tau+\copa{i}\qquad i=1,\dots,m+1,\\
b^{m+j}&=\tau_j\qquad j=2,3.
\end{split}
\end{equation}

Differently from $S^m\times S^1$, the 1-form $\tau$ is
not closed, so structure equations are a bit more complicated:
\begin{equation}\label{struc3}
\begin{split}
db^i&=b^i\wedge{\tau}+2x_ib^{m+2}\wedge b^{m+3}\qquad
i=1,\dots,m+1,\\
db^{m+2}&=2b^{m+3}\wedge\tau,\\
db^{m+3}&=-2b^{m+2}\wedge\tau,
\end{split}
\end{equation}
where the 1-form $\tau$ is related to ${\cal B}^*$ by
\[
\tau=\sum_{i=1}^{m+1}x_ib^i.
\]

\begin{oss}
The same argument used above for $S^m\times S^3$ can be applied to the
Hopf fibration $S^7\rightarrow \Cp{3}$ to obtain a
frame on $S^m\times S^7$. Nevertheless, formulas in this case are much more
complicated.\EndRem
\end{oss}

Proposition \ref{parm3} and the previous remark can be easily generalized:
\begin{teo}[\cite{BruPEP}]
Let $N^n$ be any parallelizable $n$-dimensional manifold. Then $S^m\times N$ is
parallelizable.
\end{teo}

\begin{D}
Let $T=T_1,T_2,\dots,T_n$ be a frame on $N$. The required
parallelization is thus given by
\begin{equation*}
\begin{split}
b_i&\ug M_i+x_i T\qquad i=1,\dots,m+1,\\
b_{m+j}&\ug T_j\qquad j=2,\dots,n,
\end{split}
\end{equation*}
where $M_i$ is the $\imath^{\text{th}}$ meridian vector field on
$S^m$.
\end{D}


\section{The general problem:
when is a product of spheres parallelizable?}\label{bruintro}


The proof of the theorem of Kervaire cited in the introduction is
here sketched:

\begin{Dsketch}[ of theorem \ref{teokervaire} (Kervaire)]
\begin{enumerate}
\item Show by induction there exists an embedding of $S^{n_1}\times\dots\times
S^{n_r}$ in $\R^{n_1+\dots+n_r+1}$. This is true for $r=1$. Let
\[
\map{f=(f_1,\dots,f_{n_1+\dots+n_{r-1}+1})}{S^{n_1}\times\dots\times
S^{n_{r-1}}}{\R^{n_1+\dots+n_{r-1}+1}}
\]
be the embedding given by
the inductive hypothesis, where $f$ is chosen in such a way that $f_1\ge
0$. Let $u\in S^{n_1}\times\dots\times
S^{n_{r-1}}$, and let $(\xi_1,\dots,\xi_{{n_r}+1})\in S^{n_r}$: the
embedding $f$ is thus given by
\[
\begin{array}{cll}
S^{n_1}\times\dots\times S^{n_r} &
\stackrel{f}{\longrightarrow} & \R^{n_1+\dots+n_r+1}\\
(u,(\xi_1,\dots,\xi_{{n_r}+1})) & \longmapsto & (f_2(u),\dots,f_{n_1+\dots+n_{{r-1}+1}}(u),
\xi_1\sqrt{f_1(u)},\dots,\xi_{{n_r}+1}\sqrt{f_1(u)});
\end{array}
\]
\item
suppose without any loss of generality
that the odd dimension is not $n_1$, and observe that
the degree of the Gauss map of the embedding $f$ built in 1.\ is
given by 
\[
\chi(D^{n_1+1}\times S^{n_2}\times\dots\times S^{n_r})=
\chi(D^{n_1+1})\chi(S^{n_2})\dots\chi(S^{n_r})=0,
\]
where $D^{n_1+1}$ denotes a topological disk of dimension $n_1+1$;
\item
denote by $G_{k,n}$ and $V_{k,n}$ the Grassmannian
and the Stiefel-Whitney manifold of oriented $k$-planes  and
oriented orthonormal frames in $\R^{k+n}$,
respectively. The tangential map
\[
S^{n_1}\times\dots\times S^{n_r}\longrightarrow G_{n_1+\dots+n_r,1}
\]
is null-homotopic, since by 2.\ the Gauss map is;
\item
last, denote by $P(S^{n_1}\times\dots\times S^{n_r})$ the
principal bundle of $S^{n_1}\times\dots\times S^{n_r}$, and look
at the following diagram to end the proof:
\[
\xymatrix{
P(S^{n_1}\times\dots\times S^{n_r})\ar@{.>}[r] \ar@{.>}[d] &
V_{n_1+\dots+n_r,1}\ar[d] \\
S^{n_1}\times\dots\times S^{n_r}\ar[r] &
G_{n_1+\dots+n_r,1}
}
\]
\end{enumerate}
\end{Dsketch}



Note that, due to the
homotopy theory considerations, the above proof is not very suitable
to write down explicit parallelizations on products of spheres.


Another proof of Kervaire's theorem can be developed using
a series of hints contained in the book \cite[exercises 3,4,5 and 6 of section
4.2]{HirDTo}. Details of such a proof, as developed by the
author, are given in the following.


%
%
%

In what follows, $\varepsilon_B^k$ denotes
the trivial vector bundle of rank $k$ with base
space $B$; moreover, whenever $\al$ is a vector bundle,
$E(\al),p_{\al},B(\al)$ denote the total space, the projection and
the base space of \al\ respectively.


\begin{lem}\label{lembun}
Let $\al$ be a vector bundle.
The Whitney sum $\al\oplus\varepsilon_{B(\al)}^k$ is described by
\[
\begin{split}
E(\al\oplus\varepsilon_{B(\al)}^k)&\iso E(\al)\times\R^k,\\
p_{\al\oplus\varepsilon_{B(\al)}^k}(e,v)&=p_{\al}(e),\\
B(\al\oplus\varepsilon_{B(\al)}^k)&=B(\al).
\end{split}
\]
\end{lem}

\begin{D}
The Whitney sum $\al\oplus\varepsilon_{B(\al)}^k$ is given by the pull-back
of $\al\times\varepsilon_{B(\al)}^k$ by means of the diagonal map
$B(\al)\rightarrow B(\al)\times
B(\al)$ (see for instance \cite[page
27]{MiSCCl}). Then
\[
E(\al\oplus\varepsilon_{B(\al)}^k)
=\{(e,b,v,b)\in E(\al)\times B(\al)\times
\R^k\times B(\al)\st
p_{\al}(e)=b\}
\]
and the thesis follows.
\end{D}

\begin{cor}\label{lemmabundles}
Let $\al,\be$ be vector bundles. Then, for any $k\ge 0$,
\[
\al\times(\be\oplus\varepsilon_{B(\be)}^k)\iso
(\al\oplus\varepsilon_{B(\al)}^k)\times\be.
\]
\end{cor}

\begin{D}
Observe that
\[
\begin{split}
E(\al\times(\be\oplus\varepsilon_{B(\be)}^k))&\iso E(\al)\times
E(\be\oplus\varepsilon_{B(\be)}^k)\stackrel{\ref{lembun}}{\iso} E(\al)\times
E(\be)\times\R^k,\\
E(\al\oplus\varepsilon_{B(\al)}^k)\times\be&\iso E(\al\oplus\varepsilon_{B(\al)}^k) \times E(\be)
\stackrel{\ref{lembun}}{\iso} E(\al)\times\R^k\times
E(\be),
\end{split}
\]
and use the obvious isomorphism.
\end{D}

\begin{teo}\label{probundles}
Suppose $M^m$ and $N^n$ satisfy the following
properties:
\begin{enumerate}
\item $T(M)\oplus\varepsilon_M^1$ is trivial;
\item $T(N)\oplus\varepsilon_N^1$ is trivial;
\item there is a non-vanishing vector field on $N$.
\end{enumerate}
Then $M\times N$ is parallelizable.
\end{teo}

\begin{D}
Let $\nu$ be a complement in $T(N)$ of the non-vanishing vector field on
$N$, that is,
\begin{equation}\label{c}
T(N)\iso\nu\oplus\varepsilon_N^1.
\end{equation}
Then
\begin{equation}\label{chain}
\begin{split}
T(M\times N)&\iso T(M)\times T(N)\stackrel{\myref{c}}{\iso}
T(M)\times(\nu\oplus\varepsilon_N^1)\\
&\stackrel{\ref{lemmabundles}}{\iso}
(T(M)\oplus\varepsilon_M^1)\times\nu
\stackrel{\text{1.}}{\iso}
\varepsilon_M^{m+1}\times\nu\\
&\stackrel{\ref{lemmabundles}}{\iso}\varepsilon_M^{m-1}\times(\nu\oplus\varepsilon_N^2)
\stackrel{2.}{\iso}
\varepsilon_M^{m-1}\times\varepsilon_N^{n+1}
\end{split}
\end{equation}
\end{D}

\begin{oss}\label{ossshort}
Whenever $N$ is itself parallelizable, formula \myref{chain} can
be shortened:
\begin{equation}\label{chainshort}
\begin{split}
T(M\times N)&\iso T(M)\times T(N)\stackrel{}{\iso}
T(M)\times\varepsilon_N^{n}\\
&\stackrel{}{\iso}
(T(M)\oplus\varepsilon_M^1)\times\varepsilon_N^{n-1}
\stackrel{}{\iso}
\varepsilon_M^{m+1}\times\varepsilon_N^{n-1}.
\end{split}
\end{equation}\EndRem
\end{oss}

The embedding $S^n\subset\R^{n+1}$ gives the triviality of
$T(S^n)\oplus\varepsilon_{S^n}^1$; whenever $n$ is odd, a non-vanishing vector
field on $S^n\subset\C^{(n+1)/2}$
is given by the complex multiplication. Thus, the following:
\begin{cor}\label{corbundles}
Let $n$ be any positive odd integer. Then the manifold $S^m\times S^n$
is parallelizable.
\end{cor}

And finally:

{\noindent\nopagebreak\em Second proof of theorem \ref{teokervaire}: }
Apply $r-1$ times the corollary \ref{lemmabundles} to show
that $T(S^{n_2}\times\dots\times S^{n_r})\oplus\varepsilon_{S^{n_2}\times\dots\times
S^{n_r}}^1$ is a trivial vector bundle, and use theorem
\ref{probundles}.
\EndDim

\section{An explicit parallelization $\mathcal{P}$ for products of 2
spheres}\label{secparexp}

An explicit parallelization $\mathcal{B}$ has already been found on $S^m\times S^n$,
for $n=1,3,7$, in the previous sections. Can one use theorem
\ref{probundles} to explicitly find a parallelization on any
parallelizable $S^m\times
S^n$? Answer is positive.

The trick in theorem \ref{probundles} is simple: split $TN$ by
means of the never-vanishing vector field, then use the trivial
summand to parallelize $TM$, and last detach a rank 2 trivial
summand to parallelize the remaining part of $TN$. Remark
\ref{ossshort} simply says that if $N$ is itself parallelizable,
one can avoid to detach the rank 2 trivial summand from $M$, using
the parallelization of $N$ instead.

{\em Here and henceforth, $n$ is supposed to be the odd dimension in $S^m\times S^n$.}



Denote by $y=(y_j)$
the coordinates on $\R^{n+1}$, and let $S^n\subset\R^{n+1}$ be
given by
\[
\begin{split}
S^n&\ug\{y=(y_1,\dots,y_{n+1})\in\R^{n+1}\st
y_1^2+\dots+y_{n+1}^2=1\}.
\end{split}
\]

Being $n$ odd, a never-vanishing vector field, and hence a versor
field,
is defined on $S^n$: {\em
here and henceforth, $T$ denotes the versor field  on $S^n$
given by
multiplication by $i$ in $\C^{(n+1)/2}$, namely,
\begin{equation}\label{eqtjexp}
T\ug
-y_2\pa{y_1}+y_1\pa{y_2}+\dots-y_{n+1}\pa{y_n}+y_n\pa{y_{n+1}}.
\end{equation}
}
When a shorter form of $T$ is needed, $t_j$ denotes the
coordinates of $T$, that is,
\begin{equation}
T=
\sum_{j=1}^{n+1}t_j\pa{y_j}
\end{equation}
where $t_j$ is given by
\begin{equation}\label{eqtj}
t_j=\left\{\begin{split}
-&y_{j+1}\qquad \text{if }j\text{ is
odd,}\\
&y_{j-1}\qquad \text{if }j\text{ is
even.}
    \end{split}\right.
\end{equation}

Moreover, denote by $N$ the normal versor field of
$S^n\subset\R^{n+1}$ (recall that $M$ denotes the normal versor field of
$S^m\subset\R^{m+1}$):
\begin{equation}
\begin{split}
N\ug\sum_{j=1}^{n+1}y_j\pa{y_j}.
\end{split}
\end{equation}

It is convenient to think of $T(S^m\times
S^n)=TS^m\times TS^n$ as a Riemannian subbundle of
$T\R^{m+1}_{|_{S^m}}\times T\R^{n+1}_{|_{S^n}}$;
this last is trivial, and an orthonormal frame is
$\{\partial_{x_1},\dots,\partial_{x_{m+1}},\partial_{y_1},\dots,\partial_{y_{n+1}}\}$.

Denote by $N_j$ the $\jmath^{\text{th}}$ meridian vector field on
$S^n$ (recall that $M_i$ denotes the $\imath^{\text{th}}$ meridian vector field on
$S^m$):
\[
\begin{split}
N_j&\ug\text{orthogonal
projection of $\pa{y_j}$ on $S^{n}$}\qquad j=1,\dots,n+1.
\end{split}
\]

The tangent space in a point  $(x,y)\in S^m\times S^n$ is thus
given by an euclidean  vector subspace
\[
T_xS^m\oplus T_yS^n\subset
\R^{m+1}\oplus\R^{n+1},
\]
which is generated by the
$m+n+2$ vectors $\{M_1(x),\dots,M_{m+1}(x),N_1(y),\dots,N_{n+1}(y)\}$.

One also has
\begin{equation}\label{uga}
T_xS^m\oplus\langle M(x)\rangle_\R
=\R^{m+1}\qquad\text{and}\qquad T_yS^n\oplus\langle N(y)\rangle_\R
=\R^{n+1}.
\end{equation}

As in formula \myref{fundeq}, one obtains
\begin{equation}\label{fundeqold}
\begin{split}
\pa{x_i}&=M_i+x_iM\qquad i=1,\dots,m+1,\\
\pa{y_j}&=N_j+y_jN\qquad j=1,\dots,n+1.
\end{split}
\end{equation}

Moreover, denote by $T(y)^\perp$ the vector subspace of $T_y(S^n)$
which is orthogonal to $T(y)$:
\begin{equation}\label{ugb}
\langle T(y)\rangle_\R\oplus T(y)^\perp=T_yS^n.
\end{equation}


{\em In what follows, some
computation on the vector space $T_x(S^m)\oplus
T_y(S^n)$ is done. For the sake of simplicity, the argument of
vector fields is omitted, that is, $T$ stands for $T(y)$, $M$ stands for
$M(x)$ etc\dots}

Formula \myref{chain} in theorem \ref{probundles} gives the following chain
of pointwise isomorphisms:
\newlength{\temp}
\settowidth{\temp}{$\scriptstyle\myref{ugb}$}
\begin{equation}\label{chainlinear}
\begin{split}
T_x(S^m)\oplus T_y(S^n)&\stackrel{\myref{ugb}}{=}
T_x(S^m)\oplus\langle T\rangle_\R\oplus T^\perp\\
&\stackrel{\makebox[\temp]{$\scriptstyle\al$}}{\iso}
T_x(S^m)\oplus\langle M\rangle_\R\oplus T^\perp\\
&\stackrel{\myref{uga}}{=}
\R^{m+1}
\oplus T^\perp\\
&\stackrel{\makebox[\temp]{$\scriptstyle\be$}}{\iso}\R^{m-1}\oplus \langle N\rangle_\R\oplus\langle T\rangle_\R
\oplus
T^\perp\\
&\stackrel{\myref{ugb}}{=}\R^{m-1}\oplus \langle N\rangle_\R\oplus T_yS^n\\
&\stackrel{\myref{uga}}{=}
\R^{m-1}\oplus\R^{n+1},
\end{split}
\end{equation}
where $\al$
is defined by
\[
\al(T)\ug M,
\]
and $\be$ is defined by
\[
\be(\pa{x_m})\ug N,\qquad\be(\pa{x_{m+1}})\ug T.
\]

Pulling back to $T_x(S^m)\oplus T_y(S^n)$ the $m-1$ generators
$\{\partial_{x_1},\dots,\partial_{x_{m-1}}\}$
of $\R^{m-1}$ one obtains
\begin{equation}\label{framea}
\begin{split}
\pa{x_i}&\stackrel{\myref{fundeq}}{=}
M_i+x_iM\\
&\stackrel{\al^{-1}}{\longmapsto}M_i+x_iT
\end{split}\qquad i=1,\dots,m-1,
\end{equation}
whereas pulling back to $T_x(S^m)\oplus T_y(S^n)$ the $n+1$ generators
$\{\partial_{y_1},\dots,\partial_{y_{n+1}}\}$
of $\R^{n+1}$ one obtains the more complicated formulas
\begin{equation}\label{frameb}
\begin{split}
\pa{y_j}&\stackrel{\myref{fundeq}}{=}
N_j+y_jN\\
&\stackrel{\phantom{\myref{fundeq}}}{=}\langle N_j,T\rangle T+(N_j-\langle N_j,T\rangle T)+y_jN\\
&\stackrel{\be^{-1}}{\longmapsto}
\langle N_j,T\rangle \pa{x_{m+1}}+(N_j-\langle N_j,T\rangle T)+y_j\pa{x_{m}}\\
&\stackrel{\myref{fundeq}}{=}
\langle N_j,T\rangle (M_{m+1}+x_{m+1}M)+(N_j-\langle N_j,T\rangle T)+y_j(M_{m}+x_{m}M)\\
&\stackrel{\al^{-1}}{\longmapsto}
\langle N_j,T\rangle (M_{m+1}+x_{m+1}T)+(N_j-\langle N_j,T\rangle T)+y_j(M_{m}+x_{m}T)
\end{split}j=1,\dots,n+1.
\end{equation}

The following theorem applies the above argument to $S^m\times
S^n$, odd $n$, in order to obtain an explicit frame on it:
\begin{teo}
Let $n$ be odd, and let $T=\sum_{j=1}^{n+1}t_j\pa{y_j}$ be the
tangent versor field on $S^n$
given by formula \myref{eqtj}.
Also, let $\{M_i\}_{i=1,\dots,m+1}$ and $\{N_j\}_{j=1,\dots,n+1}$
be the meridian vector fields on $S^m$
and $S^n$ respectively. Last, let $M$ and $N$ be the normal versor
fields of $S^m\subset\R^{m+1}$
and $S^n\subset\R^{n+1}$ respectively.
The product $S^m\times S^n$ is parallelized by the
 frame ${\cal P}\ug\{p_1,\dots,p_{m+n}\}$ given by
\begin{equation}\label{frame}
\begin{split}
p_i&\ug M_i+x_iT\qquad i=1,\dots,m-1,\\
p_{m-1+j}&\ug y_j M_{m}+t_j M_{m+1}+(t_j x_{m+1}+y_j
x_{m}-t_j)T+N_j
\qquad j=1,\dots,n+1.
\end{split}
\end{equation}
Moreover, $\cal P$ is orthonormal with respect to the standard metric on $S^m\times S^n$.
\end{teo}

\begin{D}
Observe that
\[
\langle N_j,T\rangle\stackrel{\myref{fundeq}}{=}\langle \pa{y_j}-y_jN,T\rangle=\langle \pa{y_j},T\rangle=t_j\qquad j=1,\dots,n+1
\]
and use
formulas \myref{framea} and
\myref{frameb}
to obtain \myref{frame}.
The orthonormality can be proved by observing that both \al\ and
\be\ in \myref{chainlinear} are isometries. But one can also
directly
check the $p_i$'s,
taking into
account formula \myref{fundeqscal}.
\end{D}

\section{The frames $\cal P$ and $\cal B$ on $S^m\times S^1$ and  $S^m\times S^3$}

If $n=1,3$ or $7$, remark \ref{ossshort} can be used to obtain a
parallelization simpler than $\cal P$
on $S^m\times S^n$. If $n=1,3$ this parallelization is just the one given in
sections $\ref{genm1}, \ref{genm3}$ respectively, which was called $\cal B$.
In this section
relations between $\cal B$ and $\cal P$ are exploited.

Let $n=1$.
Formula \myref{frame} gives the frame ${\cal P}=\{p_1,\dots,p_{m+1}\}$
on $S^m\times S^1$, whereas the frame $\cal B$ is given by formula \myref{progeofor}.
Clearly,
\[
p_i=b_i\qquad i=1,\dots,m-1.
\]
Since $\pa{\T}=-y_2\pa{y_1}+y_1\pa{y_2}=T$, one obtains
\[
\begin{split}
\langle N_1,\pa{\T}\rangle&=\langle \pa{y_1}-y_1N,-y_2\pa{y_1}+y_1\pa{y_2}\rangle=-y_2,\\
\langle N_2,\pa{\T}\rangle&=\langle \pa{y_2}-y_2N,-y_2\pa{y_1}+y_1\pa{y_2}\rangle=y_1,
\end{split}
\]
and thus
\[
\begin{split}
N_1&=-y_2T,\\
N_2&=y_1T.
\end{split}
\]
Whence
\[
\begin{split}
p_{m}&=y_1(M_{m}+x_{m}T)-y_2(M_{m+1}+x_{m+1}T)+y_2T-y_2T=y_1b_m-y_2b_{m+1},\\
p_{m+1}&=y_2(M_{m}+x_{m}T)+y_1(M_{m+1}+x_{m+1}T)-y_1T+y_1T=y_2b_m+y_1b_{m+1},
\end{split}
\]
and one gets
\begin{equation}\label{changebtop}
{\cal P}={\cal B}\left(\begin{array}{ccc|cc}
 & & & 0 & 0\\
 &I_{m-1}& &\vdots &\vdots \\
 & &  & 0& 0\\ \hline
 0&\cdots &0 &  y_1 & y_2 \\
 0&\cdots &0  & -y_2 & y_1
                \end{array}\right)
\end{equation}


Brackets of $\cal P$ are thus easily obtained by means of formulas
\myref{changebtop}, \myref{bracket1gen} and \myref{kobano}:
\begin{equation}\label{bracket1genp}
\begin{split}
[p_i,p_j]&=x_ip_j-x_jp_i\qquad i,j=1,\dots,m-1\\
[p_i,p_m]&=(-x_my_1+x_{m+1}y_2)p_i+x_ip_m-x_ip_{m+1}\qquad
i=1,\dots,m-1\\
[p_i,p_{m+1}]&=(-x_my_2-x_{m+1}y_1)p_i+x_ip_m+x_ip_{m+1}\qquad
i=1,\dots,m-1\\
[p_m,p_{m+1}]&=(x_m(y_1-y_2)-x_{m+1}(y_1+y_2))p_m +(x_m(y_1+y_2)+x_{m+1}(y_1-y_2))p_{m+1}
\end{split}
\end{equation}

Formula \myref{changebtop} gives the frame ${\cal P}^*$
dual to $\cal P$:
\[
\begin{split}
p^i&=b^i\qquad i=1,\dots,m-1,\\
p^{m}&=y_1b^m-y_2b^{m+1},\\
p^{m+1}&=y_2b^m+y_1b^{m+1}.
\end{split}
\]

The structure equations for $\cal P$ are thus obtained
by a straightforward computation:
\begin{equation}\label{diff1genp}
\begin{split}
dp^i&=\copa{i}\wedge\tau=p^i\wedge{\tau}\qquad i=1,\dots,m-1,\\
dp^m&=p^m\wedge\tau+p^{m+1}\wedge\tau,\\
dp^{m+1}&=p^{m+1}\wedge\tau-p^m\wedge\tau,
\end{split}
\end{equation}
where
$\tau$ is given by
\[
\tau=\sum_{i=1}^{m+1}x_ib^i=\sum_{i=1}^{m-1}x_ip^i+(x_my_1-x_{m+1}y_2)p^m+
(x_my_2+x_{m+1}y_1)p^{m+1}.
\]


Let $n=3$.
Formula \myref{frame} gives the frame ${\cal P}=\{p_1,\dots,p_{m+3}\}$
on $S^m\times S^3$, whereas the frame $\cal B$ is given by formula \myref{parsms3}.
Clearly,
\[
p_i=b_i\qquad i=1,\dots,m-1.
\]
Denote by ``$(*)_{\text{$\jmath^{\text{th}}$}}$'' the
$\jmath^{\text{th}}$ coordinate of $*$.
Since
\[
\begin{split}
\langle N_j-t_jT,T\rangle&=0\\
\langle N_j-t_jT,b_{m+2}\rangle&=(b_{m+2})_{\text{$\jmath^{\text{th}}$}},\\
\langle N_j-t_jT,b_{m+3}\rangle&=(b_{m+3})_{\text{$\jmath^{\text{th}}$}},\\
\end{split}\qquad j=1,\dots,4,
\]
one gets
\[
p_{m-1+j}=y_jb_m+t_jb_{m+1}+(b_{m+2})_{\text{$\jmath^{\text{th}}$}}b_{m+2}
+(b_{m+3})_{\text{$\jmath^{\text{th}}$}}b_{m+3}\qquad j=1,\dots,4.
\]
Whence
\begin{equation}\label{changebtop3}
{\cal P}={\cal B}\left(\begin{array}{ccc|cccc}
 & & & 0 & 0 & 0 & 0\\
 &I_{m-1}& &\vdots &\vdots &\vdots &\vdots \\
 & &  & 0& 0 & 0 & 0\\ \hline
 0&\cdots &0 &  y_1 & y_2 & y_3 & y_4 \\
 0&\cdots &0  & -y_2 & y_1 & -y_4 & y_3 \\
  0&\cdots &0  & -y_3 & y_4 & y_1 & -y_2 \\
 0&\cdots &0  & -y_4 & -y_3 & y_2 & y_1 \\
                \end{array}\right)
\end{equation}


Brackets of $\cal P$ can be obtained by means of a not straightforward computation
using formulas
\myref{changebtop3}, \myref{bracket3gen} and \myref{kobano}. One can
also refer to the next section, where general formulas for
$\mathcal{P}$ are given.


Formula \myref{changebtop3} gives the frame ${\cal P}^*$
dual to $\cal P$:
\begin{equation*}
{{\cal P}^*}={{\cal B}^*}\left(\begin{array}{ccc|cccc}
 & & & 0 & 0 & 0 & 0\\
 &I_{m-1}& &\vdots &\vdots &\vdots &\vdots \\
 & &  & 0& 0 & 0 & 0\\ \hline
 0&\cdots &0 &  y_1 & y_2 & y_3 & y_4 \\
 0&\cdots &0  & -y_2 & y_1 & -y_4 & y_3 \\
  0&\cdots &0  & -y_3 & y_4 & y_1 & -y_2 \\
 0&\cdots &0  & -y_4 & -y_3 & y_2 & y_1 \\
                \end{array}\right)
\end{equation*}

\section{General formulas for $\mathcal{P}$}

Recall that $T=\sum_{j=1}^{n+1}t_j\pa{y_j}$. Set
\[
\begin{split}
X_m&\ug M_m+x_mT,\\
X_{m+1}&\ug M_{m+1}+x_{m+1}T,\\
C_{j,k}&\ug y_jt_k-y_kt_j\qquad j,k=1,\dots,n+1,\\
D_{j,k}&\ug 2C_{j,k}\underbrace{\mp\delta_{k,j\pm 1}}_{j\text{ }\frac{\text{odd}}{\text{even}}}
\underbrace{\pm\delta_{j,k\pm 1}}_{k\text{ }\frac{\text{odd}}{\text{even}}}
\qquad j,k=1,\dots,n+1.
\end{split}
\]

Formulas \myref{frame} easily give
\[
\begin{split}
\sum_{j=1}^{n+1}y_jp_{m-1+j}&=M_m+x_mT=X_m,\\
\sum_{j=1}^{n+1}t_jp_{m-1+j}&=M_{m+1}+x_{m+1}T=X_{m+1}.
\end{split}
\]
A hard calculation then gives
\begin{equation}\label{genbracket}
\begin{split}
[p_i,p_{j}]&=x_ip_j-x_jp_i\qquad i,j=1,\dots,m-1,\\
[p_i,p_{m-1+j}]&=-(y_jx_m+t_jx_{m+1})p_i\\
&\underbrace{\mp
x_i p_{m-1+j\pm 1}}_{j\text{ }\frac{\text{odd}}{\text{even}}}
+x_iy_jX_m+x_it_jX_{m+1}\qquad i=1,\dots,m-1,j=1,\dots,n+1,\\
[p_{m-1+j},p_{m-1+k}]&=D_{j,k}\sum_{i=1}^{m-1}x_ip_i+y_jp_{m-1+k}-y_kp_{m-1+j}\\
&+(x_mD_{j,k}-x_{m+1}C_{j,k})X_m+((x_{m+1}-1)D_{j,k}+x_mC_{j,k})X_{m+1}\\
&+\underbrace{(\mp y_jx_m\mp t_jx_{m+1}\pm
t_j)p_{m-1+k\pm 1}}_{k\text{ }\frac{\text{odd}}{\text{even}}}\\
&+\underbrace{(\pm y_kx_m\pm t_kx_{m+1}\mp
t_k)p_{m-1+j\pm 1}}_{j\text{ }\frac{\text{odd}}{\text{even}}}
\qquad j,k=1,\dots,n+1.\\
\end{split}
\end{equation}


\begin{acknowledgements}
This paper is part of a Ph.~D.\
thesis, and the author wishes in a special way
to thank Paolo Piccinni for the
motivation and the constant help.
\end{acknowledgements}



\end{document}